%---------------------------------------------------------
% wild3.tex
%
%                   Ringel-Schmidmeier
%
%   Date: February,21, 2004  September 5, 2004  To be translated using AMSTeX
%
%   ------------------------------------------------------
%        Submodule Categories of Wild Representation Type
%   ------------------------------------------------------
%
%---------------------------------------------------------
\magnification=\magstep1   
\input pictex
\input amstex
\UseAMSsymbols
\hoffset=0truecm 
\hsize=125mm 
\vsize=165mm
\NoBlackBoxes
\nopagenumbers
\parindent=8mm
\mathsurround=1pt
\vcorrection{1cm}
\font\gross=cmbx10 scaled\magstep1   \font\abs=cmcsc10
\font\rmk=cmr8  \font\itk=cmti8    \font\ttk=cmtt8
%=========================

 \def \End{\operatorname{End}}
 \def \rad{\operatorname{rad}}

 \def \sub{\subseteq} 
 
 \def \Hom{\operatorname{Hom}}

 \def \dim{\operatorname{dim}}

 \def \mod{\operatorname{mod}}

 \def\arr#1#2{\arrow <2mm> [0.25,0.75] from #1 to #2}
 \def\sq{\plot 0 0  1 0  1 1  0 1  0 0 /}
 \def\halfsq{\plot 0 0  .5 0  .5 1  0 1  0 0 /}

          %Fuer Text innerhalb von $$...$$. Verwendung: \T{oder}
 \def\E#1{{\parindent=1truecm\narrower\narrower\noindent{#1}\smallskip}}  
          %Zum Einruecken, Verwendung: \E{dieses wird eingerueckt}
 \def\qed{\phantom{m.} $\!\!\!\!\!\!\!\!$\nolinebreak\hfill\checkmark}

%=================================
\headline{\ifnum\pageno=1\hfill %
    \else\ifodd\pageno \Rechts \else \Links \fi  \fi}
    \def\Links{\abs \the\pageno\hfil Ringel, Schmidmeier\hfil}
    \def\Rechts{\abs \hfil Wild Representation Type\hfil\the\pageno}
%============================================
%               Version    Datum
%                  !         !
%                  v         v
%
%\noindent{\rmk [wild3, Sept 5, 2004]}
%
%============================================
%
        \vglue1truecm
\centerline{\gross Submodule Categories of Wild Representation Type}
        \bigskip
\centerline{Claus Michael Ringel and Markus Schmidmeier}
        \bigskip
%============================================

\E{{\it Abstract.}
Let $\Lambda$ be a commutative local uniserial ring 
of length at least seven  with radical factor ring $k$.
We consider the category $\Cal S(\Lambda)$ of all
possible embeddings of submodules of finitely generated $\Lambda$-modules
and show that $\Cal S(\Lambda)$ 
is controlled $k$-wild with a single control object $I\in\Cal S(\Lambda)$. 
In particular, it follows that each finite dimensional $k$-algebra 
can be realized as a quotient
$\End(X)/\End(X)_I$ 
of the endomorphism ring of some object $X\in\Cal S(\Lambda)$
modulo the ideal $\End(X)_I$ of all maps which factor through a finite direct
sum of copies of $I$.}

\bigskip\medskip

Let $\Lambda$ be a ring. 
Recall that an object $M = (M_0,M_1)$ (or written also $(M_1\subseteq M_0)$)
in the submodule category 
$\Cal S(\Lambda)$ consists of a finitely generated $\Lambda$-module $M_0$
together with a $\Lambda$-submodule $M_1$ of $M_0$; a morphism $f\:M \to N$
in $\Cal S(\Lambda)$ 
is just a $\Lambda$-linear map $f\:M_0\to N_0$  which
preserves the submodules, that is, $f(M_1)\subseteq N_1$ holds. 

In this paper, $\Lambda$ always will be a commutative local uniserial ring
of finite length $n$. Usually, we will assume that $n \ge 7$.
The radical factor field will be denoted by $k$ 
and $t$ will be a radical generator (thus $\Lambda/\langle t\rangle = k$).
We have the following two special cases in mind: First of all, 
if $\Lambda$ is the ring $\Bbb Z/\langle p^n\rangle$ where $p$ is a prime number, then we 
are dealing with the category of all possible embeddings of a subgroup in 
a $p^n$-bounded finite abelian group; the classification problem for the
objects in $\Cal S(\Bbb Z/\langle p^n\rangle )$ was raised by Birkhoff [B] in 1934.
Second, if $\Lambda = k[T]/\langle T^n\rangle$, 
where $k[T]$ is the polynomial ring in
one variable $T$ over the field $k$, then we consider the possible invariant
subspaces of a nilpotent operator (indeed, the objects in $\Cal 
S(k[T]/\langle T^n\rangle)$
may be written as triples $(V,\phi,U)$, where $V$ is a k-space, $\phi\:V \to V$
is a $k$-linear transformation with $\phi^n = 0$ und $U$ is a subspace of $V$
with $\phi(U) \subseteq U$). 

Some remarks concerning notions of ``wildness'' of additive categories
will be given in the last sections. In the case $\Lambda = 
\Bbb Z/\langle p^n\rangle$, 
Arnold [A] has shown that $\Cal S(\Bbb Z/\langle p^{10}\rangle)$ is ``wild''.
In the case $\Lambda = k[T]/\langle T^n\rangle $ Simson has shown in [Si] that 
$\Cal S(k[T]/\langle T^7\rangle )$ is ``wild'' whereas 
$\Cal S(k[T]/\langle T^6\rangle )$ is still tame, thus providing the precise bound for ``wildness''.
It is not surprising that the special case $\Lambda = k[T]/\langle T^n\rangle$ is better 
understood, since in this case many powerful techniques are available
(in particular covering theory). The main result presented here will not 
dependent on $\Lambda$ being an algebra over a field.
In particular, it applies to the 
classical case of subgroups of finite abelian groups, as considered by Birkhoff, and it can
be used in order to construct parametrized families of metabelian groups [Sc].
In case $\Lambda$ is an algebra over a field, the last section shows in which way
the main result can be strengthened.
      \bigskip\medskip

\centerline{\bf Controlled Wildness}
        \medskip
Let $\Cal A$ be an additive category and $\Cal C$ a class of objects (or
a full subcategory) in $\Cal A$. 
Given objects $A,A'$ in $\Cal A$, we will write 
$\Hom(A,A')_{\Cal C}$ for the set of maps $A \to A'$ which factor through a 
(finite) direct sum of objects in $\Cal C$ (note that in this way we 
attach to $\Cal C$ the ideal $\langle\Cal C\rangle$
in $\Cal A$ generated by the identity
morphisms of the objects in $\Cal C$). 
The same convention will apply to a single object
$C$ in $\Cal A$: we denote by $\Hom(A,A')_C$ the set of 
maps $A \to A'$ which factor through a 
(finite) direct sum of copies of $C$.
Of course, given an ideal $\Cal I$ of $\Cal A$, we write $\Cal A/\Cal I$
for the corresponding factor category: it has the same objects as
$\Cal A$ and given two objects $A,A'$ of $\Cal A$, the group
$\Hom_{\Cal A/\Cal I}(A,A')$ is defined as $\Hom_{\Cal A}(A,A')/\Cal I(A,A').$
In particular, the category $\Cal A/\langle\Cal C\rangle$ has the same objects as $\Cal A$
and 
$$
 \Hom_{\Cal A/\langle\Cal C\rangle}(A,A') =\Hom_{\Cal A}(A,A')/\Hom(A,A')_{\Cal C}.
$$
        \medskip\noindent
{\bf Definition.\/}
We say that $\Cal A$ is {\it controlled $k$-wild} provided there are
full  subcategories $\Cal C \subseteq \Cal B \subseteq \Cal A$ 
such that $\Cal B/\langle\Cal C\rangle$ is equivalent to 
$\mod k\langle X,Y\rangle$ where
$k\langle X,Y\rangle$ is the free $k$-algebra in
two generators. 
We will call $\Cal C$ the {\it control class}, 
and in case $\Cal C$ is given by
a single object $C$ then this object $C$ will be the {\it control
object.} We refer to [R] for a discussion of controlled wildness.
        \medskip\bigskip
\centerline{\bf The Setting}
        \medskip
We are going to show that the category
$\Cal S(\Lambda)$ is controlled $k$-wild. In order to do so, we
need to find suitable full subcategories $\Cal C \subseteq \Cal B
\subseteq \Cal S(\Lambda)$. In fact, $\Cal C$ will consists of a single
object $I$, whereas $\Cal B$ will be a suitable subcategory of
the ``interval'' in-between
the object $I$ and a related one $J$ with $I\subset J$.
Given two objects $I \subset J$ in $\Cal S(\Lambda)$, we denote by
the {\it interval} $[I,J]$ the class of all objects $M$ of
$\Cal S(\Lambda)$ such that $I^m \subseteq M \subseteq J^m$
for some natural number $m$.
        \medskip
In order to exhibit objects in $\Cal S(\Lambda)$, it is convenient
to use some graphical description. It is well-known and easy to see that
the indecomposable $\Lambda$-modules are up to isomorphism of the form 
$\Lambda/\langle t^i\rangle$ with
$1 \le i \le n$, thus the indecomposable $\Lambda$-modules are
characterized by the length ($\Lambda/\langle t^i\rangle $ has length $i$). The
Krull-Remak-Schmidt theorem asserts that the isomorphism classes of the
$\Lambda$-modules (of finite length) correspond bijectively to the
partitions $\lambda = (\lambda_1,\dots,\lambda_m)$ with all parts
$\lambda_i \le n$; the $\Lambda$-module corresponding to the partition
$(\lambda_1,\dots,\lambda_m)$ is just 
$\bigoplus_i \Lambda/\langle t^{\lambda_i}\rangle$, or,
equivalently, the $\Lambda$-module
with generators $x_1,\dots,x_m$ and defining relations $t^{\lambda_i}x_i = 0$,
for $1\le i \le m$. 
We will attach to a partition its Young diagram using an
arrangement of boxes, however we will deviate from the usual convention as
follows: the various parts will be drawn vertically and not horizontally,
and the parts will not necessarily be adjusted at the top or the socle.
For example, we will consider below the partition (7,4,2), and it will be
suitable to draw the corresponding Young diagram as follows:
$$
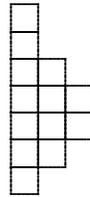
\hbox
{\beginpicture 
\setcoordinatesystem units <0.3cm,0.3cm>
\multiput{\sq} at 0 6 0 5  0 4  0 3  0 2  0 1  0 0  1 4  1 3  1 2  1 1  2 3  2 2  /
\endpicture}
$$      
The left column has 7 boxes, the middle one 4 and the right column 2 boxes,
as the partition $(7,4,2)$ asserts.
The adjustment of these columns made here
depends on the fact that we have in mind a particular submodule, 
and we want that there is a generating system for the submodule such that any of
these
generators is a linear combination of elements which belong to boxes at the same
height.
        \medskip
Here are the objects $I$ and $J$: The $\Lambda$-module
$J_0$ is given by the partition $(7,4,2)$, say with generators $x,y,z$,
annihilated by $t^7, t^4, t^2$ respectively, and $I_0$ is generated by
$tx,y,z$, thus it corresponds to the partition $(6,4,2).$ The submodule $J_1$
is generated by $t^3x-ty$ and $ty-z$, and $I_1 = tJ_1$. 
$$
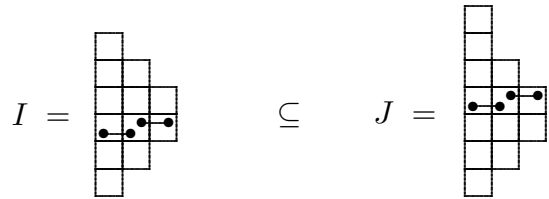
\hbox{\beginpicture 
\setcoordinatesystem units <0.3cm,0.3cm>
\put{$I\;=$} at -2 2.5
\multiput{\sq} at 0 5  0 4  0 3  0 2  0 1  0 0  1 4  1 3  1 2  1 1  2 3  2 2  /
\put{$\ssize \bullet$} at 0.3 1.8 
\put{$\ssize \bullet$} at 1.3 1.8 
\put{$\ssize \bullet$} at 1.7 2.2 
\put{$\ssize \bullet$} at 2.7 2.2 
\plot 0.3 1.8  1.3 1.8 /
\plot 1.7 2.2  2.7 2.2 /
\endpicture}\qquad\quad
\hbox{\beginpicture 
\setcoordinatesystem units <0.3cm,0.3cm>
\put{$\subseteq \qquad J\;=$} at -4 2.5
\multiput{\sq} at 0 6 0 5  0 4  0 3  0 2  0 1  0 0  1 4  1 3  1 2  1 1  2 3  2 2  /
\put{$\ssize \bullet$} at 0.3 2.8 
\put{$\ssize \bullet$} at 1.3 2.8 
\put{$\ssize \bullet$} at 1.7 3.2 
\put{$\ssize \bullet$} at 2.7 3.2 
\plot 0.3 2.8  1.3 2.8 /
\plot 1.7 3.2  2.7 3.2 /
\endpicture}
$$      
In these pictures, we have indicated the generators of the submodules 
using pairs of bullets which are connected by a horizontal line (and the 
shift of the
columns was accomplished  in such a way that the connecting lines become
horizontal lines).
        \medskip
We are going to describe the interval $[I,J]$ in-between
the objects $I$ and $J$ in terms of representations of a quiver 
$\Delta$. The quiver $\Delta$ looks as follows:
$$
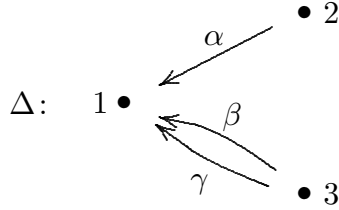
\hbox{\beginpicture\setcoordinatesystem units <2mm,2mm>
\put{$\Delta\:$} at -5 0
\put{$\bullet$} at 0 0
\put{$\bullet$} at 10 5
\put{$\bullet$} at 10 -5
\put{$1$} at -1.3 0
\put{$2$} at 11.5 5
\put{$3$} at 11.5 -5
\put{$\alpha$} at 5 3.6 
\put{$\beta$}  at 6 -.8
\put{$\gamma$} at 4.2 -4.5
\arr{8.2 4.2} {2 1}
\arr{3.5 -2.4} {1.8 -1.2}
\arr{4 -1.2} {2 -0.8}
\setquadratic
\plot 8 -4.6  5  -3.3  3.5 -2.4 /
\plot 8.4 -3.7  6  -2.1  4 -1.2 /
\endpicture}
$$
It has three vertices: one sink (labelled $1$) and two sources
(labelled $2$ and $3$), thus there are three simple representations
$S(1), S(2), S(3)$. The simple representation $S(1)$ is projective, the simple
representations $S(2)$ and $S(3)$ are injective. 
Let us denote by $\mod_e k\Delta$ the full subcategory of $\mod k\Delta$
given by all representations without a simple direct summand. Note that the representations of $\Delta$ without a simple injective
direct summand are precisely the socle-projective representations 
(of course, a
representation is said to be {\it socle-projective} provided the socle is
projective). We denote by $\mod_{\text{sp}} k\Delta$ the full subcategory of all
socle-projective representations. The inclusion functors
$$
 \mod_{e}k\Delta \subset 
 \mod_{\text{sp}}k\Delta \subset 
 \mod k\Delta 
$$
allow us to identify the categories
$$
 \mod_{e}k\Delta =  
 \mod_{\text{sp}}k\Delta /\langle S(1)\rangle  =
 \mod k\Delta /\langle S(1),S(2),S(3)\rangle,
$$
since all
the simple representations of $\Delta$ are projective or injective.
Note that $\mod k\Delta / \langle S(1),S(2),S(3)\rangle$ is the factor category
of $\mod k\Delta$ modulo the ideal of maps which factor through semisimple objects. 
         \medskip
The key to proving the controlled wildness of $\Cal S(\Lambda)$ is the following result.

\medskip\noindent
{\bf Theorem 1.} {\it Let $\Lambda$ be a
commutative local uniserial ring
of length $n \ge 7$ and let $k$ be its
radical factor field. Then, the factor category 
$[I,J]/\langle I\rangle$ is
equivalent to the category $\mod_e k\Delta$.}
        \bigskip

The definition of $[I,J]$ may be rephrased as follows: In order to form the
direct sums $I^m, J^m$ of copies of $I$ and $J$, respectively, 
let $\widetilde W$
be a free $\Lambda$-module of rank $m$ so that we may identify the inclusion
$I^m \subset J^m$ with the map 
$\widetilde W\otimes I \to \widetilde W\otimes J$ 
induced by the inclusion $I \subset J$. 
Then each object $M$ in $[I,J]$ can be visualized as the middle term 
in a sequence of inclusions of the following type:
$$
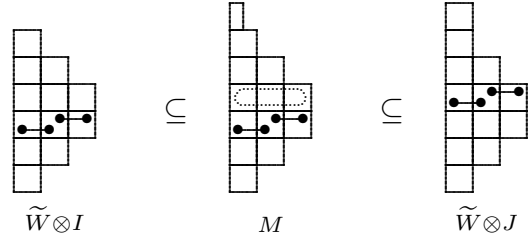
\hbox{\beginpicture 
\setcoordinatesystem units <0.3cm,0.3cm>
\put{} at -2 2.5
\multiput{\sq} at 0 5  0 4  0 3  0 2  0 1  0 0  1 4  1 3  1 2  1 1  2 3  2 2  /
\put{$\ssize \bullet$} at 0.3 1.8 
\put{$\ssize \bullet$} at 1.3 1.8 
\put{$\ssize \bullet$} at 1.7 2.2 
\put{$\ssize \bullet$} at 2.7 2.2 
\plot 0.3 1.8  1.3 1.8 /
\plot 1.7 2.2  2.7 2.2 /
\put{$\ssize \widetilde W\otimes I$} at 1.5 -1.5
\endpicture}\qquad
\hbox{\beginpicture 
\setcoordinatesystem units <0.3cm,0.3cm>
\put{$\subseteq$} at -2 2.5
\multiput{\sq} at 0 5  0 4  0 3  0 2  0 1  0 0  1 4  1 3  1 2  1 1  2 3  2 2  /
\put{\halfsq} at 0 6
\put{$\ssize \bullet$} at 0.3 1.8 
\put{$\ssize \bullet$} at 1.3 1.8 
\put{$\ssize \bullet$} at 1.7 2.2 
\put{$\ssize \bullet$} at 2.7 2.2 
\plot 0.3 1.8  1.3 1.8 /
\plot 1.8 2.2  2.6 2.2 /
\setdots<1.3pt>
\plot .2 2.9  .2 3.1  .4 3.3  2.6 3.3  2.8 3.1  2.8 2.9  2.6 2.7  .4 2.7  .2 2.9 /
\put{$\ssize M$} at 1.5 -1.7
\endpicture}\qquad
\hbox{\beginpicture 
\setcoordinatesystem units <0.3cm,0.3cm>
\put{$\subseteq$} at -2 2.5
\multiput{\sq} at 0 6 0 5  0 4  0 3  0 2  0 1  0 0  1 4  1 3  1 2  1 1  2 3  2 2  /
\put{$\ssize \bullet$} at 0.3 2.8 
\put{$\ssize \bullet$} at 1.3 2.8 
\put{$\ssize \bullet$} at 1.7 3.2 
\put{$\ssize \bullet$} at 2.7 3.2 
\plot 0.3 2.8  1.3 2.8 /
\plot 1.7 3.2  2.7 3.2 /
\put{$\ssize \widetilde W\otimes J$} at 1.5 -1.5
\endpicture}
$$
Here, the dotted region represents the quotient 
$M_1/(\widetilde W\otimes_\Lambda I_1)$; and the half box on the top 
corresponds to the quotient $M_0/(\widetilde W\otimes_\Lambda I_0)$.
(Note that now the boxes no longer correspond to individual composition factors
of the $\Lambda$-module $M_0$, but to suitable semisimple subfactors.)
        \bigskip\medskip
\centerline{\bf The Layer Functors}

\nopagebreak\medskip
Let us analyse the objects $M$ in $[I,J]$.
There are the {\it layer functors} 
$$
 L_i\: [I,J] \to \mod \Lambda
$$ 
defined by
$$
\align
 L_1M &= t^4M_0 \cap t^{-1}0 \cr
 L_2M &= t^3M_0 \cap t^{-2}0 \cr
 L_iM &= t^{-i+2}L_{2} \qquad\qquad\text{for } i \ge 3.
\endalign
$$ 
Note that the $L_iM$ are $\Lambda$-submodules of $M_0$. They
form a filtration of $M_0$ and  
can be visualized as follows:
$$
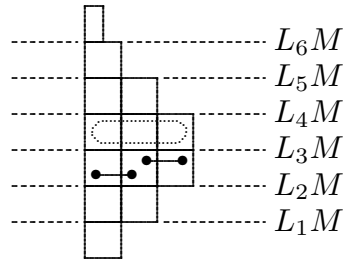
\hbox
{\beginpicture 
\setcoordinatesystem units <0.4cm,0.4cm>
\put{} at 0 0
\multiput{\sq} at 0 5  0 4  0 3  0 2  0 1  0 0  1 4  1 3  1 2  1 1  2 3  2 2  /
\put{\halfsq} at 0 6
\put{$\ssize \bullet$} at 0.3 1.8 
\put{$\ssize \bullet$} at 1.3 1.8 
\put{$\ssize \bullet$} at 1.7 2.2 
\put{$\ssize \bullet$} at 2.7 2.2 
\plot 0.3 1.8  1.3 1.8 /
\plot 1.8 2.2  2.6 2.2 /
\setdots<1.3pt>
\plot .2 2.9  .2 3.1  .4 3.3  2.6 3.3  2.8 3.1  2.8 2.9  2.6 2.7  .4 2.7  .2 2.9 /
\setdashes<.5mm>
\plot -2 5.5  -0.2 5.5 /
\plot -2 4.5  -0.2 4.5 /
\plot -2 3.5  -0.2 3.5 /
\plot -2 2.5  -0.2 2.5 /
\plot -2 1.5  -0.2 1.5 /
\plot -2 0.5  -0.2 0.5 /
\plot  1.2 5.5  5 5.5 /
\plot  2.2 4.5  5 4.5 /
\plot  3.2 3.5  5 3.5 /
\plot  3.2 2.5  5 2.5 /
\plot  3.2 1.5  5 1.5 /
\plot  2.2 0.5  5 0.5 /

\put{$L_6M$} at 6.2 5.5
\put{$L_5M$} at 6.2 4.5
\put{$L_4M$} at 6.2 3.5
\put{$L_3M$} at 6.2 2.5
\put{$L_2M$} at 6.2 1.5
\put{$L_1M$} at 6.2 0.5
\endpicture}
$$
This definition of the submodules $L_iM$ 
only depends on $M_0$, it does not take into account $M_1$. 

\smallskip Of special
interest is the following observation:

\medskip\noindent{\bf Lemma.} 
{\it The subobject $(M_1\cap L_3M\subseteq L_6M)$ of $M$
is a direct
sum of copies of $I$, and any homomorphism from $I$ to 
$M$ maps into $(M_1\cap L_3M\subseteq L_6M)$.} 

\smallskip We call this subobject  
$(M_1\cap L_3M\subseteq L_6M)$ the {\it $I$-socle} of $M$. 
        \smallskip\noindent
{\it Proof:\/} Let $\widetilde W$ be a free $\Lambda$-module such that
the inclusions $\widetilde W \otimes I\sub M\sub \widetilde W\otimes J$ hold.
The inclusion $\widetilde W\otimes I \sub M$
embeds $\widetilde W\otimes I$ into
$(M_1\cap L_3M\subseteq L_6M)$ and clearly $\widetilde W\otimes I_0 = L_6M.$
But we also have $\widetilde W \otimes I_1 = \widetilde W\otimes tJ_1 =
M_1\cap L_3M.$ This shows that   $\widetilde W\otimes I = 
(M_1\cap L_3M\subseteq L_6M)$, thus 
$(M_1\cap L_3M\subseteq L_6M)$ is a direct sum of copies of $I$.
Given a map $I \to M,$ it will send $I_0 = L_6I$ into $L_6M$ and
$I_1 = I_1\cap L_3I$ into $M_1\cap L_3M$, thus it maps into
$(M_1\cap L_3M\subseteq L_6M)$. \qed

\bigskip\medskip
\centerline{\bf From $\Cal S(\Lambda)$ to Representations of $\Delta$}
        \medskip
Given $M$ in $[I,J]$, let $F(M)$ be defined by
$$
\hbox{\beginpicture\setcoordinatesystem units <2mm,2mm>
\put{$L_1M$} at -2 0
\put{$M/L_6M$} at 11 5
\put{$M_1/(M_1\cap L_3M)$} at 15 -5
\put{$\alpha$} at 5 3.6 
\put{$\beta$}  at 6 -.8
\put{$\gamma$} at 4.2 -4.5
\arr{8 4} {2 1}
\arr{3.5 -2.4} {1.8 -1.2}
\arr{4 -1.2} {2 -0.8}
\setquadratic
\plot 8 -4.6  5  -3.3  3.5 -2.4 /
\plot 8.4 -3.7  6  -2.1  4 -1.2 /
\endpicture}
$$
with $\alpha = t^6,$  $\beta = t^3,$ 
and a $k$-linear map $\gamma$ which still has to be specified.

Actually, let us define a surjective homomorphism $\gamma'\:L_4M \to L_1M$
with kernel $L_3M + pL_5M$, the required map $\gamma$ will be induced
by the restriction of $\gamma'$ to $M_1$ (note that $M_1 \subseteq 
L_4M$). The map $\gamma'$ will yield an isomorphism between 
the following two
shaded boxes:
$$
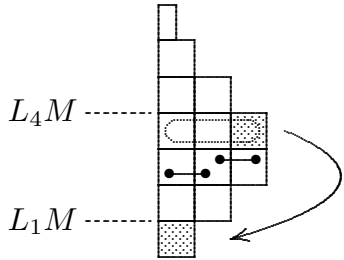
\hbox
{\beginpicture 
\setcoordinatesystem units <0.4cm,0.4cm>
\put{} at 0 0
\multiput{\sq} at 0 5  0 4  0 3  0 2  0 1  0 0  1 4  1 3  1 2  1 1  2 3  2 2  /
\put{\halfsq} at 0 6
\put{$\ssize \bullet$} at 0.3 1.8 
\put{$\ssize \bullet$} at 1.3 1.8 
\put{$\ssize \bullet$} at 1.7 2.2 
\put{$\ssize \bullet$} at 2.7 2.2 
\plot 0.3 1.8  1.3 1.8 /
\plot 1.8 2.2  2.6 2.2 /

\setdots<1pt>
\plot .2 2.9  .2 3.1  .4 3.3  2.6 3.3  2.8 3.1  2.8 2.9  2.6 2.7  .4 2.7  .2 2.9 /

\setshadegrid span <.4mm>
\vshade 0 -.5 .5  1 -.5 .5 /
\vshade 2 2.5 3.5 3 2.5 3.5 /
\setsolid
\setquadratic
\plot 3.5 3  5 1.5  2.1 0.025 /
\arr{2.1 0.025}{2 0} 
\setlinear
\setdashes<.5mm>

\plot -2 3.5  -0.2 3.5 /
\plot -2 0.5  -0.2 0.5 /

\put{$L_4M$} at -3.2 3.5
\put{$L_1M$} at -3.2 0.5
\endpicture}
$$
Here is the definition of $\gamma'(c)$ for $c\in L_4M$ in a condensed form:
$$
 \gamma'(c) \quad=\quad t^2\;\bigg(\Big(\big((tc+t^2L_5M)\cap t^{-1}0\big)
+\big(M_1\cap L_3M\big)\Big)\cap t^3L_6M\bigg)
$$
(note that $\gamma'$ depends only on $L_6$ and $M_1\cap L_3M$, thus on the $I$-socle of
$M$). 
        \medskip
In order to understand the definition and to see that $\gamma'$ is really
a $\Lambda$-homomorphism, we proceed stepwise:
Thus, we start with $c\in L_4M$. Take an element $c'\in (tc+t^2L_5M)\cap t^{-1}0$,
such an element exists since $tL_4M = t^2L_5M+t^{-1}0.$
Next, take an element
$c'' \in (c'+(M_1\cap L_3M))\cap t^3L_6M$ --- again, we note that such an element 
exists, now we use that $c'$ belongs to $L_3M$ and that
$L_3M = (M_1\cap L_3M)+ t^3L_6M.$ The proposed definition of $\gamma'(c)$
amounts to $\gamma'(c) = t^2c''$. Now $c'$ is unique up the addition of elements
from $t^2L_5M\cap t^{-1}0 = t^{-1}0 \cap L_2M$, thus $c''$ is unique up to the
addition of elements from $(M_1\cap L_3M)\cap t^3L_6M = L_2M$, and the latter
elements go to zero under the multiplication by $t^2$. This shows that
$\gamma'(c)$ is a well-defined element and since $c''$ belongs to $L_3M$,
we see that $\gamma'(c)$ belongs to $L_1M$.
Of course, it is clear that such a construction yields a
homomorphism $\gamma'$. One finally verifies that $\gamma'$
is surjective and that its kernel is 
$L_3M + tL_5M$.

\medskip It also follows from the construction that a 
homomorphism $g\:M\to N$
in $\Cal S(\Lambda)$ between objects $M,N\in[I,J]$
commutes with $\gamma'$.  Hence we obtain a functor
$F\:[I,J]\to \mod k\Delta$.

\medskip\noindent{\it Example.} Under this functor $F$, 
the object $I$ is sent to $F(I) = S(1)$,
whereas $F(J)$ is the injective envelope of $S(1)$:
$$
\hbox{\beginpicture\setcoordinatesystem units <1.5mm,1.5mm>
\put{$F(I) =$} at -10 0
\put{$k$} at -2 0
\put{$0$} at 11 5
\put{$0$} at 11 -5
\arr{8 4} {2 1}
\arr{3.5 -2.4} {1.8 -1.2}
\arr{4 -1.2} {2 -0.8}
\setquadratic
\plot 8 -4.6  5  -3.3  3.5 -2.4 /
\plot 8.4 -3.7  6  -2.1  4 -1.2 /
\endpicture}
\qquad\qquad
\hbox{\beginpicture\setcoordinatesystem units <1.5mm,1.5mm>
\put{$F(J)=$} at -10 0
\put{$k$} at -2 0
\put{$k$} at 11 5
\put{$k\oplus k$} at 11 -5
\put{$1$} at 5 3.6 
\put{$\big({1\atop0}\big)$}  at 6.5 -.1
\put{$\big({0\atop1}\big)$} at 3.7 -5
\arr{8 4} {2 1}
\arr{3.5 -2.4} {1.8 -1.2}
\arr{4 -1.2} {2 -0.8}
\setquadratic
\plot 8 -4.6  5  -3.3  3.5 -2.4 /
\plot 8.4 -3.7  6  -2.1  4 -1.2 /
\endpicture}
$$
        \medskip
\noindent{\bf Proposition 1.} {\it The functor $F$ is a full 
and dense functor
from $[I,J]$ onto the category $\mod_{\text{sp}}k\Delta$ of
socle-projective representations of
$\Delta$. The representation 
$F(I) = S(1)$ is simple projective, and the kernel of the
induced functor $[I,J] \to \mod_{\text{sp}} k\Delta/\langle S(1)\rangle$
is just the ideal of all maps which factor through a
direct sum of copies of $I$.}

\smallskip 
The proof of Proposition 1 will be given at the end of the next section; 
Theorem 1 is an immediate consequence of Proposition 1.

\bigskip\medskip
\centerline{{\bf $\ldots$ and Back to $\Cal S(\Lambda)$}}

        \nopagebreak\medskip
In order to show that the functor $F$ is full and dense, we are going 
to present an inverse construction which we label $\Phi$. 
We work in the homomorphism category
$\Cal H(\Lambda)$ for $\Lambda$. The
objects in $\Cal H(\Lambda)$ are the $\Lambda$-linear maps, say 
$A = (A_1\overset a \to \to A_0)$,
and a morphism between two such objects $A = (A_1\overset a \to \to A_0)$ 
and $B = (B_1\overset b \to \to B_0)$
consists of two homomorphisms $f_0\:A_0\to B_0$ and $f_1\:A_1\to B_1$ 
such that $f_0a = bf_1$ holds.  
Clearly, $\Cal S(\Lambda)$ is just the full 
exact subcategory
of $\Cal H(\Lambda)$ of those objects $A = (A_1\overset a \to \to A_0)$ for which
the map $a$ is monic. 
        \medskip
Note that the inclusion
$I\to J$ gives rise to the short exact sequence in $\Cal H(\Lambda)$
$$
 \varepsilon\: \quad 0\longrightarrow I\longrightarrow J\longrightarrow
 (k\oplus k\overset 0 \to {\to} k)\longrightarrow 0.
$$
        \medskip
Let $W$ be a vector space, $V$ a subspace of $W$, $U$ a subspace of
$W\oplus W$, then we may consider the tripel $(W,V,U)$
as a representation of the quiver $\Delta$ as follows
$$
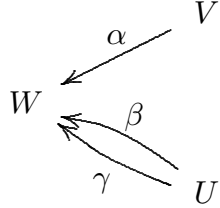
\hbox{\beginpicture\setcoordinatesystem units <2mm,2mm>
\put{$W$} at 0 0
\put{$V$} at 10 5
\put{$U$} at 10 -5
\put{$\alpha$} at 5 3.6 
\put{$\beta$}  at 6 -.8
\put{$\gamma$} at 4.2 -4.5
\arr{8 4} {2 1}
\arr{3.5 -2.4} {1.8 -1.2}
\arr{4 -1.2} {2 -0.8}
\setquadratic
\plot 8 -4.6  5  -3.3  3.5 -2.4 /
\plot 8.4 -3.7  6  -2.1  4 -1.2 /
\endpicture}
$$
with $\alpha$ the inclusion map, and $\beta$ the first,  $\gamma$ the second
projection of $U$ into $W$ (thus $\beta(w_1,w_2) = w_1, 
\gamma(w_1,w_2) = w_2$,
where $w_1,w_2 \in W$ and $(w_1,w_2)\in U$). Note that in this way,
we obtain precisely all the representations of $\Delta$ which do not
have a simple injective direct summand. 
        \medskip

Let $\widetilde W$ be a free $\Lambda$-module with $\widetilde W/\rad \widetilde W = W$. 
In the category $\Cal H(\Lambda)$, we consider
the following fibre product construction of $\widetilde W\otimes_\Lambda \varepsilon$
along the inclusion $(U\overset 0 \to {\to} V)\longrightarrow
(W\oplus W\overset 0 \to {\to} W)$:
$$\CD 0 @>>> \widetilde W\otimes_\Lambda I @>>> \Phi(W,V,U) @>>> 
                             (U\overset 0 \to {\to} V) @>>> 0\cr
   @.          @|                  @VVV              @VVV          @.   \cr
0     @>>> \widetilde W\otimes_\Lambda I @>>> \widetilde W\otimes_\Lambda J @>>> 
                            (W\oplus W \overset 0 \to {\to} W) @>>> 0
\endCD
$$
In this way, we define the object $\Phi(W,V,U)$. Note that
by the five lemma, the vertical map in the center of the above diagram
is monic, so $\Phi(W,V,U)$, being a subobject of an object in 
$\Cal S(\Lambda)$, also lies in $\Cal S(\Lambda)$ and clearly $F\Phi(W,V,U)$
is the subobject $(W,V,U)$ of $F(\widetilde W\otimes J)=(W,W,W\oplus W)$.

Let us look again at the visualization of the objects in $[I,J]$
considered above:
$$
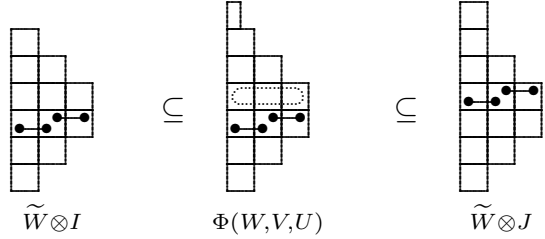
\hbox{\beginpicture 
\setcoordinatesystem units <0.3cm,0.3cm>
\put{} at -2 2.5
\multiput{\sq} at 0 5  0 4  0 3  0 2  0 1  0 0  1 4  1 3  1 2  1 1  2 3  2 2  /
\put{$\ssize \bullet$} at 0.3 1.8 
\put{$\ssize \bullet$} at 1.3 1.8 
\put{$\ssize \bullet$} at 1.7 2.2 
\put{$\ssize \bullet$} at 2.7 2.2 
\plot 0.3 1.8  1.3 1.8 /
\plot 1.7 2.2  2.7 2.2 /
\put{$\ssize \widetilde W\otimes I$} at 1.5 -1.5
\endpicture}\qquad
\hbox{\beginpicture 
\setcoordinatesystem units <0.3cm,0.3cm>
\put{$\subseteq$} at -2 2.5
\multiput{\sq} at 0 5  0 4  0 3  0 2  0 1  0 0  1 4  1 3  1 2  1 1  2 3  2 2  /
\put{\halfsq} at 0 6
\put{$\ssize \bullet$} at 0.3 1.8 
\put{$\ssize \bullet$} at 1.3 1.8 
\put{$\ssize \bullet$} at 1.7 2.2 
\put{$\ssize \bullet$} at 2.7 2.2 
\plot 0.3 1.8  1.3 1.8 /
\plot 1.8 2.2  2.6 2.2 /
\setdots<1.3pt>
\plot .2 2.9  .2 3.1  .4 3.3  2.6 3.3  2.8 3.1  2.8 2.9  2.6 2.7  .4 2.7  .2 2.9 /
\put{$\ssize \Phi(W,V,U)$} at 1.5 -1.7
\endpicture}\qquad
\hbox{\beginpicture 
\setcoordinatesystem units <0.3cm,0.3cm>
\put{$\subseteq$} at -2 2.5
\multiput{\sq} at 0 6 0 5  0 4  0 3  0 2  0 1  0 0  1 4  1 3  1 2  1 1  2 3  2 2  /
\put{$\ssize \bullet$} at 0.3 2.8 
\put{$\ssize \bullet$} at 1.3 2.8 
\put{$\ssize \bullet$} at 1.7 3.2 
\put{$\ssize \bullet$} at 2.7 3.2 
\plot 0.3 2.8  1.3 2.8 /
\plot 1.7 3.2  2.7 3.2 /
\put{$\ssize \widetilde W\otimes J$} at 1.5 -1.5
\endpicture}
$$
For such an object $\Phi(W,V,U)$, 
the dotted region represents $U$ which is the quotient
of $\Phi(W,V,U)_1$ modulo $\widetilde W\otimes_\Lambda I_1$; and the 
half box on the top corresponds to the 
subspace $V$ of $W$, which is the quotient of 
$\Phi(W,V,U)_0$ modulo $\widetilde W\otimes_\Lambda I_0$.

        \medskip
Suppose $(W,V,U)$ and $(W',V',U')$ are two such triples
(thus representations
of the quiver $\Delta$ without simple injective direct summands). 
Let $\widetilde W$ and $\widetilde W'$ be free 
$\Lambda$-modules with $\widetilde W/\rad \widetilde W = W$ and $\widetilde W'/\rad \widetilde W' = W'$,
respectively. A morphism 
$(W,V,U)\to (W',V',U')$ in the category $\mod \Delta$
is given by a map $g\: W\to W'$ such that $g(V)\subseteq V'$ and 
$(g\oplus g)(U) \subseteq U'.$

Such a map $g$ gives rise to a morphism 
$\Phi(g)\: \Phi(W,V,U)\to \Phi(W',V',U')$ in the category $\Cal S(\Lambda)$
which makes the following diagram commutative:
$$\hbox{\beginpicture\setcoordinatesystem units <1.1cm,0.8cm>
%==========================================
\put{$\ssize 0$} at 0.2 3
\put{$\ssize 0$} at 0.2 1
\put{$\ssize 0$} at 1.2 0
\put{$\ssize 0$} at 1.2 2
\put{$\ssize \widetilde W\otimes I$} at 2 3
\put{$\ssize \widetilde W\otimes I$} at 2 1
\put{$\ssize \widetilde W'\otimes I$} at 3 2
\put{$\ssize \widetilde W'\otimes I$} at 3 0
\put{$\ssize \widetilde W\otimes J$} at 4 1
\put{$\ssize \Phi(W,V,U)$} at 4 3
\put{$\ssize \widetilde W'\otimes J$} at 5 0
\put{$\ssize \Phi(W',V',U')$} at 5 2
\put{$\ssize (U\to V)$} at 6 3
\put{$\ssize (W\oplus W\to W)$} at 6 1
\put{$\ssize (U'\to V')$} at 7 2
\put{$\ssize (W'\oplus W'\to W')$} at 7 0
\put{$\ssize 0$} at 7.8 3
\put{$\ssize 0$} at 7.8 1
\put{$\ssize 0$} at 8.8 2
\put{$\ssize 0$} at 8.8 0
\arr{2.3 2.7} {2.7 2.3}
\arr{2.3 0.7} {2.7 0.3}
\arr{4.3 2.7} {4.7 2.3}
\arr{4.3 0.7} {4.7 0.3}
\arr{6.3 2.7} {6.7 2.3}
\arr{6.3 0.7} {6.7 0.3}
\arr{0.7 3} {1.3 3}
\arr{0.7 1} {1.3 1}
\arr{1.7 2} {2.3 2}
\arr{1.7 0} {2.3 0}
\arr{2.7 3} {3.1 3}
\arr{3.1 1} {3.3 1}
\plot 2.7 1  2.9 1 /
\arr{3.7 2} {4.1 2}
\arr{3.7 0} {4.3 0}
\arr{4.9 3} {5.3 3}
\arr{4.6 1} {4.9 1}
\arr{5.9 2} {6.3 2}
\arr{5.7 0} {6.1 0}
\arr{6.7 3} {7.3 3}
\arr{7.1 1} {7.3 1}
\plot 6.8 1  6.9 1 /
\arr{7.7 2} {8.3 2}
\arr{7.9 0} {8.3 0}
\arr{2 1.9} {2 1.5}
\arr{3 1.5} {3 0.5}
\arr{4 1.9} {4 1.5}
\arr{5 1.5} {5 0.5}
\arr{6 1.9} {6 1.5}
\arr{7 1.5} {7 0.5}
\plot 2 2.5  2 2.1 /
\plot 4 2.5  4 2.1 /
\plot 6 2.5  6 2.1 /
\endpicture}
$$
Indeed, one uses the projectivity of $\widetilde W$ as a $\Lambda$-module 
in order to obtain a lifting $\widetilde g$ of $g$
which makes the diagram 
$$
 \CD 
 \widetilde W                @>\widetilde g>>      \widetilde W'         \cr 
 @V{\text{can}}VV              @VV{\text{can}}V \cr 
 W                @>>g>                 W'
 \endCD
$$
commutative, then the bottom part and the right hand square in the 
three dimensional diagram commute.  Now a fibre product construction in
the category of morphisms in $\Cal H(\Lambda)$ gives the forward pointing
map in the middle of the top part of the diagram.
This is the map $\Phi(g)$ we are looking for.  
Since the vertical maps in the middle are both monic,
$\Phi(g)\:\Phi(W,V,U)\to \Phi(W',V',U')$ is just the restriction of the 
map $\widetilde W\otimes_\Lambda J\overset{\widetilde g\otimes 1}\to{\longrightarrow} 
\widetilde W'\otimes_\Lambda J$ and clearly $F\Phi(g) = g.$
        \medskip
Let us stress that the construction $\Phi$ is not functorial, since
it depends on a choice of liftings: We had to write the vector space $W$
as $W = \widetilde W/\rad \widetilde W$ for some free $\Lambda$-module $\widetilde W$ and given the linear
transformation $g\:W \to W'$, we used a lifting $\widetilde g\: \widetilde W \to \widetilde W'$
of $g$. 
        \medskip
\noindent
{\it Proof of Proposition 1:\/}
Note that $F(I) = (k,0,0)$ is the simple projective representation
of $\Delta$. Therefore, $F$ induces a functor 
$[I,J]/\langle I\rangle$ into $\mod_{\text{sp}} k\Delta/\langle S(1)\rangle$ 
and this functor $[I,J]/\langle I\rangle \longrightarrow \mod_{\text{sp}} k\Delta/
\langle S(1)\rangle$ is full and dense. 
It remains to determine its kernel. For this, let $M,M'$ be in
$[I,J]$ and consider a map $f\:M \to M'$ such that $F(f)$
factors through a direct sum of copies of $S(1)$.
It follows
that $F(f)_2 = F(f)_3 = 0.$ 
Now $F(M) = (L_1M,M/L_6M,M_1/(M_1\cap L_3M))$ and 
$F(M') = (L_1M',M'/L_6M',M'_1/(M'_1\cap L_3M'))$. The maps 
$F(f)_2\:M/L_6M \to M'/L_6M'$ and
$F(f)_3\:M_1/(M_1\cap L_3M) \to M'_1/(M'_1\cap L_3M')$ are induced
by $f$; since these are zero maps, we see that
$$
 f(M_0) \subseteq L_6M' \qquad\text{and}\qquad f(M_1) \subseteq M_1'\cap L_3M',
$$
thus $f$ maps into the $I$-socle of $M'$.
By the Lemma, the $I$-socle of $M'$ is a direct sum of copies of $I$, thus
$f$ belongs to $\Hom(M,M')_I$, as we wanted to show. \qed

        \bigskip\medskip
\centerline{\bf Conclusion}
        \medskip
Our main result is 
        \medskip\noindent
{\bf Theorem 2.} {\it
Let $\Lambda$ be a commutative local uniserial 
ring of length $n\ge 7$ and let $k$ be its 
radical factor. Then the category $\Cal S(\Lambda)$ is
controlled $k$-wild.}

        \bigskip
For the proof, we need the (well-known) fact that the category
$\mod k\Delta$ is strictly $k$-wild:
recall that an additive category $\Cal A$ is said to be {\it strictly $k$-wild} provided
there exists a full embedding of the category $\mod k\langle X,Y\rangle$
into $\Cal A$.
The following embedding $G$ of $\mod k\langle X,Y\rangle$ int
the category of representations of $\Delta$ is known to be full and exact:
consider
a $k\langle X,Y\rangle$-module $(V;X,Y)$ (here, $V$ is a $k$-space, and
$X$ and $Y$ are linear transformations of $V$, they are given by 
the multiplication using the corresponding generators with
the same names); under $G$ we sent it 
to the following representation of $\Delta$
$$
\hbox{\beginpicture\setcoordinatesystem units <2mm,2mm>
\put{$V\oplus V$} at -1 0
\put{$V$} at 10 5
\put{$V\oplus V$} at 12 -5
\put{$\alpha$} at 5 3.6 
\put{$\beta$}  at 6 -.8
\put{$\gamma$} at 4.2 -4.5
\arr{8 4} {2 1}
\arr{3.5 -2.4} {1.8 -1.2}
\arr{4 -1.2} {2 -0.8}
\setquadratic
\plot 8 -4.6  5  -3.3  3.5 -2.4 /
\plot 8.4 -3.7  6  -2.1  4 -1.2 /
\put{with} at 18 0
\put{$\alpha = \bmatrix 1 \cr 0 \endbmatrix,\ 
  \beta = \bmatrix 1 & 0 \cr 0 & 1 \endbmatrix,\ 
  \gamma = \bmatrix 0 & X \cr 1 & Y \endbmatrix$.} at 37 0
\endpicture}
$$
Note that {\it no representation in the image of $G$ has a simple direct
summand.} 
        \bigskip
        \noindent
{\it Proof of Theorem 2:\/} 
Let $\Cal B$ be the full subcategory of all objects $M$
in $[I,J]$ such that either $M = I$ or else $F(M)$ lies in the
image of the functor $G\:\mod k\langle X,Y\rangle \to
\mod_{\text{sp}} k\Delta$. 
The required equivalence $\Cal B/\langle I \rangle \to
\mod k\langle X,Y\rangle$ 
is given by the 
restriction of the
functor $F$ in Proposition 1 to $\Cal B$. \qed
        \bigskip
Let us mention some more details of this equivalence
$$
 \Cal B/\langle I \rangle \longrightarrow \mod k\langle X,Y\rangle.
$$
The $k\langle X,Y\rangle$-module $(V;X,Y)$ corresponds to
$\Phi(V\oplus V,V\oplus 0,U_{XY})$ in $[I,J]$, where
$$
 U_{XY} = \{(v_1,v_2,Xv_2,v_1+Yv_2)\mid v_1,v_2\in V\}
  \subseteq V \oplus V \oplus V \oplus V;
$$ 
the $\Lambda$-module $\Phi(V\oplus V,V\oplus 0,U_{XY})_0$
is given by the partition $(7^d,6^d,4^{2d},2^{2d})$ with  $d=\dim V_k$, 
and its
submodule $\Phi(V\oplus V,V\oplus 0,U_{XY})_1$
is given by the partition $(4^{2d},2^{2d}).$
        \bigskip
The above equivalence has the following consequence:
        \medskip
\noindent
{\bf Corollary.} {\it Let $R$ be a finite-dimensional $k$-algebra.
There exists $M$ in $\Cal S(\Lambda)$ such that $\End(M)/\End(M)_I$
is isomorphic to $R$.} 

        \bigskip
On the other hand, we stress the following (clearly also well-known) fact:
        \medskip
\smallskip\noindent{\bf Proposition 2.} 
{\it The category $\Cal S(\Lambda)$ is not strictly $K$-wild,
for any field $K$.}
        
\smallskip\noindent
{\it Proof:\/} 
Assume $\Cal S(\Lambda)$ is strictly $K$-wild, for some field $K$.
There are infinitely 
many isomorphism classes of 
finite length $K\langle X,Y\rangle$-modules $M$ with 
endomorphism ring $K$, and there are pairs $M,M'$ of such modules with
$\Hom(M,M') = 0 = \Hom(M',M)$; for example, just take for $M$ and $M'$ two
non-isomorpic one-dimensional representations.  
Thus, a full embedding of $\mod K\langle X,Y\rangle$
into $\Cal S(\Lambda)$ yields an object $(A\subseteq B)$ in $\Cal S(\Lambda)$
with endomorphism ring $K\times K$. Note that the multiplication with the
radical generator $t$ of $\Lambda$ gives a nilpotent endomorphism of 
any object $(A\subseteq B)$. Thus, if $\End(A\subseteq B) = K\times K,$
then $t$ has to act by zero on $B$. However, there are 
only two indecomposables $(A \subseteq B)$ in $\Cal S(\Lambda)$ such that $t$ 
acts as zero on $B$, namely
$$S_1=(0\subseteq k)\quad\text{and}\quad S_2=(k\subseteq k).$$
As there are nonzero maps from $S_1$ to $S_2$, it follows that $K\times K$
cannot be realized as an endomorphism ring. \qed

\bigskip
\centerline{\bf $K$-Algebras.}
        \medskip
If $\Lambda$ is a $K$-algebra (for any field $K$,
not necessarily isomorphic to the radical factor of $\Lambda$) 
then the assignments
$\widetilde W=W\otimes_K\Lambda$ and $\widetilde g=g\otimes_K\Lambda$
make $\Phi\:\mod K\Delta\to \Cal S(\Lambda)$ into a functor.

\medskip\noindent{\bf Proposition 3.} {\it Assume that $\Lambda$
is a $K$-algebra. 

\smallskip\item{1.} The functor $\Phi$ is exact (and additive)
and hence naturally 
equivalent to the tensor functor
$-\otimes_{K\Delta}\Phi(K\Delta)$.

\item{2.}The composition $F\circ \Phi$ is naturally 
equivalent to the identity functor on $\mod K\Delta$, and hence
$\Phi$ preserves indecomposables and reflects isomorphisms.

\item{3.} The exact embedding $\mod K\langle X,Y\rangle\to 
        \mod_{\text{sp}} K\Delta \to \Cal S(\Lambda)$, 
        makes the category $\Cal S(\Lambda)$ $K$-wild in the sense of Drozd.
\qed

}

        \bigskip\bigskip
\noindent
{\bf Acknowledgement}
\nopagebreak\medskip
One of the authors (MS) would like to thank 
Manfred Dugas (Baylor University, Texas)
for helpful discussions.  In fact, Manfred has pointed out to him that the 
usual concepts for wildness fail in the context of subgroup categories,
and this advice has motivated his research.

        \bigskip\bigskip
%-----------------------------------------------------------------------------
\frenchspacing

\noindent
{\bf References} \nopagebreak       {\baselineskip=9pt \rmk
\parindent=1.5truecm
                                         \medskip\smallskip\noindent
\item{[A]} D.\ M.\ Arnold: 
{\itk Abelian Groups and Representations of Finite Partially Ordered Sets}, 
Springer CMS Books in Mathematics (2000).

\smallskip\noindent
\item{[B]} G.~Birkhoff, {\itk Subgroups of abelian groups},
Proc.~Lond.~Math.~Soc., II. Ser. 38, 1934, 385--401.

\smallskip\noindent
\item{[R]} C.~M.~Ringel, 
{\itk Combinatorial representation theory. History and future} 
in: Representations of Algebras, Vol. I,
Proc.\ Conf.\ ICRA IX, Bejing 2000; Bejing Normal University Press, 2002,
122--144.

\smallskip\noindent
\item{[Sc]} M.~Schmidmeier, {\itk A construction of metabelian groups,}
To appear.

\smallskip\noindent
\item{[Si]} D.\ Simson, {\itk Chain categories of modules 
and subprojective representations of posets over uniserial algebras},
Rocky Mountain J.\ Math.~32, 2002, 1627--1650.
\par}
        \bigskip\medskip\noindent
{\rmk Claus Michael Ringel,
Fakult\"at f\"ur Mathematik, Universit\"at Bielefeld,
\par\noindent  POBox 100\,131, \ D-33\,501 Bielefeld 
\par\noindent {\ttk ringel\@mathematik.uni-bielefeld.de}
        \smallskip\noindent
Markus Schmidmeier, 
Department of Mathematical Sciences, Florida Atlantic University,
\par\noindent Boca Raton, Florida 33431-0991
\par\noindent {\ttk markus\@math.fau.edu}
}

\bye